\documentclass[preprint,12pt,authoryear]{elsarticle}

\usepackage[utf8]{inputenc}
\usepackage{amsmath,amsthm}
\usepackage{amssymb}
\usepackage[english]{babel}
\usepackage{enumitem}
\usepackage{mathtools}
\usepackage{verbatim}
\usepackage{parskip}
\usepackage{float}
\usepackage{tikz}
\usetikzlibrary{decorations.markings}
\usepackage{adjustbox}
\usepackage{ytableau}
\usetikzlibrary{matrix, calc, arrows}
\usepackage{comment}
\usepackage{fullpage}
\usepackage[bookmarksopen=false,pdftex=true,breaklinks=true,%
      backref=page,pagebackref=true,plainpages=false,%
      hyperindex=true,pdfstartview=FitH,colorlinks=true,%
      pdfpagelabels=true,colorlinks=true,linkcolor=blue,%
      citecolor=red,urlcolor=green,hypertexnames=false%
      ]%
   {hyperref}

\tikzstyle{vertexwhite}=[draw,thick,fill=white,circle,inner sep=2pt]
\tikzstyle{vertexblack}=[draw,thick,fill=black,circle,inner sep=2pt]

\usepackage{graphicx}

\usepackage[top=0.85in,left=0.85in,footskip=0.75in]{geometry}

\usepackage{changepage}

\usepackage{textcomp,marvosym}

\usepackage[right]{lineno}

\usepackage[nopatch=eqnum]{microtype}
\DisableLigatures[f]{encoding = *, family = * }

\newtheorem{theorem}{Theorem}[section]

\newtheorem{lemma}[theorem]{Lemma}
\newtheorem{definition}{Definition}
\newtheorem{conjecture}{Conjecture}

\newtheorem{observation}{Observation}
\begin{document}

\begin{frontmatter}

\title{Dense, irregular, yet always graphic $3$-uniform hypergraph degree sequences}

\author[1,4]{Runze Li}
\author[1,2,3]{Istv\'an Mikl\'os}

%\affil[$3$]{Budapest Semesters in Mathematics\\ 1071 Budapest, Bethlen G. t\'er 2\\ Hungary}
\affiliation[1]{organization={Budapest Semesters in Mathematics},%Department and Organization
            addressline={Bethlen G. ter 2}, 
            city={Budapest},
            postcode={1071}, 
            state={},
            country={Hungary}}

\affiliation[2]{organization={HUN-REN Renyi Institute},%Department and Organization
            addressline={Realtanoda u. 13-15}, 
            city={Budapest},
            postcode={1053}, 
            state={},
            country={Hungary}}

%\affil[$2$]{SZTAKI\\ 1111 Budapest, L\'agym\'anyosi u. 11\\ Hungary}
\affiliation[3]{organization={HUN-REN SZTAKI},%Department and Organization
            addressline={Lagymanyosi u. 11}, 
            city={Budapest},
            postcode={1111}, 
            state={},
            country={Hungary}}

%\affiliation[4]{University of California Santa Barbara, Santa Barbara, CA, U.S.}
\affiliation[4]{organization={University of California Santa Barbara},%Department and Organization
            city={Santa Barbara},
            postcode={93106}, 
            state={CA},
            country={USA}}

\begin{abstract}
    A $3$-uniform hypergraph is a generalization of simple graphs where each hyperedge is a subset of vertices of size $3$. The degree of a vertex in a hypergraph is the number of hyperedges incident with it. The degree sequence of a hypergraph is the sequence of the degrees of its vertices. The degree sequence problem for $3$-uniform hypergraphs is to decide if a $3$-uniform hypergraph exists with a prescribed degree sequence. Such a hypergraph is called a realization. Recently, Deza \emph{et al.} proved that the degree sequence problem for $3$-uniform hypergraphs is NP-complete. Some special cases are easy; however, polynomial algorithms have been known so far only for some very restricted degree sequences. The main result of our research is the following. If all degrees are between $\frac{2n^2}{63}+O(n)$ and $\frac{5n^2}{63}-O(n)$ in a degree sequence $D$, further, the number of vertices is at least $45$, and the degree sum can be divided by $3$, then $D$ has a $3$-uniform hypergraph realization. Our proof is constructive and in fact, it constructs a hypergraph realization in polynomial time for any degree sequence satisfying the properties mentioned above. To our knowledge, this is the first polynomial running time algorithm to construct a $3$-uniform hypergraph realization of a highly irregular and dense degree sequence.
\end{abstract}


\end{frontmatter}

%\maketitle
\newcommand{\5}{\left\lfloor\frac{5n^2}{7}\right\rfloor}
\newcommand{\2}{\left\lceil\frac{2n^2}{7}\right\rceil}

\section{Introduction}
A degree sequence $D = d_1, d_2, \ldots, d_n$ is a sequence of non-negative integers.
The degree sequence problem asks if there exists a simple graph $G = (V, E)$ with prescribed degrees of its vertices $D = d_1, d_2, \ldots, d_n$. The graph $G$ is called the realization of $D$.
The degree sequence problem is one of the first solved problems in algorithmic graph theory. In 1955 and in 1962, Havel \cite{Havel} 
and Hakimi 
\cite{Hakimi} 
independently gave polynomial running time algorithms to decide if a realization of a degree sequence $D$ exists and if so, the algorithm also constructs a realization. The running time of these algorithms grows polynomically with the length of the degree sequence.
A few years later, Erd{\H o}s and Gallai \cite{ErdosGallai} gave the inequalities that are necessary and sufficient to have a simple graph with a prescribed degree sequence. Gale \cite{Gale} and Ryser \cite{Ryser} gave the inequalities that are necessary and sufficient to have a bipartite graph with the prescribed degree sequences of the two vertex classes. The Havel-Hakimi algorithm \cite{Havel,Hakimi} for bipartite graphs is a folklore. 

Hypergraphs are generalizations of simple graphs. A hyperedge $e\in E$ of a hypergraph $H = (V, E)$ is a non-empty subset of $V$. A hypergraph is $k$-uniform if each edge is a subset of vertices of size $k$. 
A hyperedge  $e$ is incident with a vertex $v$ if $v\in e$, and the degree of a vertex is the number of its incident hyperedges.
The degree sequence problem can be naturally generalized to hypergraphs.
It was an open problem for a long time if polynomial running time algorithms could solve the hypergraph degree sequence problem. A few years ago, Deza \emph{et al.} \cite{Dezaetal2018,Dezaetal2019} proved that it is already NP-complete to decide if a $3$-uniform hypergraph exists with a prescribed degree sequence. On the other hand, efficient algorithms have been developed for some special classes of degree sequences. These efficient algorithms can decide if a hypergraph realization exists when the degree sequences are very close to the regular degree sequences \cite{Frosinietal,fpr,Palmaetal,aaff} or the degree sequence is sparse \cite{aldosarigreenhill}.

Since the general degree sequence problem is hard for $3$-uniform hypergraphs, it is a natural attempt to characterize the degree sequences for which the degree sequence problem can be solved in polynomial time. In algorithmic graph theory, many algorithmic problems can be solved more easily for sparse and regular degree sequences. For example, the edge packing problem is NP-hard in general \cite{durretal}. However, tree degree sequences are sparse, and it is easy to pack two trees \cite{kundu-tree}, and there are partial results to pack many tree degree sequences \cite{ghm,msww,kundu3tree}. Also, packing half-regular degree sequences is easy \cite{AMZ}. Asymptotic formulae only for the number of regular graphs \cite{bollobas} and the number of sparse $r$-uniform hypergraphs \cite{aldosarigreenhill} with a given degree sequence exist.

An intensively studied class of degree sequences are the P-stable degree sequences which have certain properties that other degree sequences do not have \cite{egmmss}. Dense degree sequences might be P-stable. For example, there exists a continuous range of $(c_1,c_2)$ pairs such that any degree sequence whose degrees are between $c_1n$ and $c_2n$ are P-stable \cite{egmmss}. We call these degree sequences \emph{linearly bounded degree sequences}. An analogous degree sequence class for $3$-uniform hypergraphs would contain degree sequences with degrees between $c_1n^2$ and $c_2n^2$. Indeed, while the degree of each vertex in a complete graph $K_n$ is $n-1$,
the degree of each vertex in the $3$-uniform complete hypergraph on $n$ vertices is ${n-1\choose 2}$. Therefore, a hypergraph degree sequence of length $n$ could be considered dense if the average degree is $\Omega(n^2)$, or even more if each degree is $\Omega(n^2)$. Similarly, a $3$-uniform hypergraph degree sequence can be considered highly irregular if the difference between its maximal and minimal degrees is $\Omega(n^2)$. In this paper, we characterize a highly irregular, dense degree sequence class that is always graphic. In other words, each degree sequence in this class has a $3$ - uniform hypergraph realization. We simply require that the minimum degree be at least $\frac{2n^2}{63}+\frac{4n}{5}+5$, the maximum degree be at most $\frac{5n^2}{63}-\frac{20n}{63}$\footnote{In Theorem~\ref{theo:hypergraph-linear-bound}, we give a slightly wider, though less readable bound}, the number of vertices be at least $45$, and the degree sum is divisible by $3$. This divisibility property is trivially necessary and can be viewed as the generalization of the well-known ``handshaking lemma" for simple graph degree sequences. Our proof is constructive and yields a polynomial running time algorithm that constructs a $3$-uniform hypergraph realization for any degree sequence in the class.

The main ingredients of the algorithm are the following. First, we introduce the tripartite, $3$-uniform hypergraphs which are generalizations of the bipartite graphs. In Section~\ref{sec:tripartite} we characterize an always graphic tripartite, $3$-uniform hypergraph degree sequence class on $n+n+n$ vertices. Particularly, we show that any tripartite $3$-uniform degree sequence is graphic if the minimum degree is at least $\frac{2n^2}{7}$, the maximum degree is at most $\frac{5n^2}{7}$, and the degree sums in the three vertex classes are the same. Then in Section~\ref{sec:linear-hypergraphic}, we show that using a few hyperedges that take only $O(n)$ degrees on each vertex and all the degrees on at most $2$ vertices, any $3$-uniform hypergraph degree sequence on $n$ vertices with the above-mentioned constraints can be tailored into a graphic tripartite, $3$-uniform hypergraph degree sequence on $\left\lfloor\frac{n}{3}\right\rfloor+\left\lfloor\frac{n}{3}\right\rfloor+\left\lfloor\frac{n}{3}\right\rfloor$ vertices.

\section{Preliminaries}
\begin{definition}
A hypergraph $H = (V, E)$ is a generalization of simple graphs. For all $e\in E$, $e$ is a non-empty subset of $V$. A hyperedge $e$ is  \emph{incident} with $v$ if $v\in e$.
A hypergraph is \emph{$t$-uniform} if for all $e\in E$, $e\in {V\choose t}$. A hypergraph $H = (V, E)$ is \emph{partite $t$-uniform} if $V$ is a disjoint partition of $V_1, V_2, \ldots, V_t$, and for all $e \in E$ and all $i=1, 2,\ldots t$, $|e\cap V_i| =1$, i.e. each edge is incident with exactly one vertex in each vertex class. 
\end{definition}
To make it simple, we call partite, $3$-uniform hypergraphs \emph{tripartite hypergraphs}. 
In this paper, any hypergraph will be either a $3$-uniform hypergraph or a tripartite hypergraph. We will omit the attributive ``$3$-uniform" when we talk about $3$-uniform hypergraphs.
\begin{definition}
The \emph{degree}  of a vertex of a hypergraph is the number of hyperedges incident with it. The \emph{degree sequence} of a hypergraph is the sequence of the degrees of its vertices. If a hypergraph is tripartite, then the degree sequence can be naturally broken down by the vertex classes, that is, it can be written as
$$
(d_{1,1}, d_{1,2}, \ldots, d_{1,n_1}), (d_{2,1}, d_{2,2},\ldots, d_{2,n_2}), (d_{3,1}, d_{3,2}, \ldots, d_{3,n_3}).
$$

A degree sequence is \emph{$k$-regular} if each degree is $k$. A degree sequence is \emph{almost regular} if each degree is either $k$ or $k+1$ for some $k$.

Given $D$ is a sequence of non-negative integers, we say that a hypergraph $H = (V, E)$ is a \emph{realization} of $D$, if the sequence of the degrees of the vertices of $H$ is $D$. If $D$ has a realization, then we say that $D$ is \emph{graphic}.
\end{definition}

Simple graphs can be considered as $2$-uniform hypergraphs and bipartite graphs can be considered as partite, $2$-uniform hypergraphs. Therefore the definitions above on degree sequences, their graphicality, and their realizations can be naturally extended to simple and bipartite graphs.

The following Lemma~\ref{lem:balancing} gives the foundation for why degree sequences with certain bounds are always graphic (given that the sum of the degrees satisfies some basic rules). Before we state and prove the lemma, we introduce the concept called hinge flip \cite{ak21}. We introduce two hinge flip operations, one on the realization of degree sequences and another on degree sequences. 
A hinge flip operation on a realization $G = (V, E)$ of a degree sequence removes a(n) (hyper)edge $\{v_i\}\cup X \in E$ and adds a(n) (hyper)edge $\{v_j\}\cup X \notin E$. The set $X$ has cardinality $2$ in the case of hypergraphs or tripartite hypergraphs and has cardinality $1$ in the case of simple or bipartite graphs.
The corresponding hinge flip operation on a degree sequence $D = (d_1, \dots d_n)$ is an operation where we decrease a $d_i$ by $1$ and increase a $d_j$ of $D$ by $1$.
It is easy to see that the modified realization will be a realization of the degree sequence modified by the corresponding hinge flip operation. Any hinge flip on a realization of a degree sequence has a corresponding hinge flip operation on the corresponding degree sequence. However, there might be a degree sequence $D$, its realization $G$, and a hinge flip operation $HF$ on $D$ such that $HF$ does not have a corresponding hinge flip on $G$.

If $d_i > d_j$, we call it \emph{balancing hinge flip}, otherwise, we call it \emph{reverse hinge flip}. Lemma~\ref{lem:balancing} states that any balancing hinge flip on a graphic degree sequence does not change the graphicality because in any realization of the degree sequence, there is a corresponding hinge flip operation.

\begin{lemma}\label{lem:balancing}
    \begin{enumerate}
        \item Let $D = (D_1,D_2)$ be a graphic bipartite degree sequence, and let $d_i,d_j \in D_1$ and $d_i<d_j$. Let $D'$ be the bipartite degree sequence obtained from $D$ by adding $1$ to $d_i$ and subtracting $1$ from $d_j$. Then any realization of $D$ has a balancing hinge flip operation yielding a realization of $D'$, thus $D'$ is also a graphic bipartite degree sequence.
      \item Let $D = (D_1,D_2,D_3)$ be a graphic tripartite hypergraph degree sequence, and let $d_i,d_j \in D_1$ and $d_i<d_j$. Let $D'$ be the tripartite hypergraph degree sequence obtained from $D$ by adding $1$ to $d_i$ and subtracting $1$ from $d_j$. Then any realization of $D$ has a balancing hinge flip operation yielding a realization of $D'$, thus $D'$ is also a graphic tripartite hypergraph degree sequence.
      \item Let $D$ be a  graphic hypergraph degree sequence, and let $d_i,d_j \in D_1$ and $d_i<d_j$. Let $D'$ be the hypergraph degree sequence obtained from $D$ by adding $1$ to $d_i$ and subtracting $1$ from $d_j$. Then any realization of $D$ has a balancing hinge flip operation yielding a realization of $D'$, thus $D'$ is also a graphic hypergraph degree sequence.
    \end{enumerate}
\end{lemma}
\begin{proof}
    \begin{enumerate}
        \item Let $G = (U,V,E)$ be a realization of $D$, and let $v_i$ and $v_j$ be the vertices whose degrees are $d_i$ and $d_j$. Since $d_j>d_i$, there is a $u$ such that $(v_j,u)\in E$ and $(v_i,u)\notin E$. Then $G' = (U,V,(E\cup \{(v_i,u)\})\setminus \{(v_j,u)\})$ is a realization of $D'$, thus $D'$ is graphic.
        \item Let $H = (A,B,C,E)$ be a realization of $D$, and let $a_i$ and $a_j$ be the vertices whose degrees are $d_i$ and $d_j$. Since $d_j>d_i$, there are $b\in B$ and $c\in C$  such that $(a_j,b,c)\in E$ and $(a_i,b,c)\notin E$. Then $H' = (A,B,C,(E\cup \{(a_i,b,c)\})\setminus \{(a_j,b,c)\})$ is a realization of $D'$, thus $D'$ is graphic.
        \item Let $H = (V,E)$ be a realization of $D$, and let $v_i$ and $v_j$ be the vertices whose degrees are $d_i$ and $d_j$. Since $d_j>d_i$, there are $v_1,v_2\in V$ such that $\{v_1,v_2\}\cap \{d_i,d_j\} = \emptyset$, and further $(v_j,v_1,v_2)\in E$ and $(v_i,v_1,v_2)\notin E$. Then $H' = (V,(E\cup \{(v_i,v_2,v_3)\})\setminus \{(v_j,v_2,v_3)\})$ is a realization of $D'$, thus $D'$ is graphic.
    \end{enumerate}    
\end{proof}
The corollary of Lemma~\ref{lem:balancing} is that any degree sequence between linear bounds is graphic if some extreme degree sequence is graphic. This is stated precisely in the following theorem.
\begin{theorem}\label{theo:balancing}
    \begin{enumerate}
        \item 
        \begin{sloppypar}
        Let $D = (D_1,D_2)$ be a graphic bipartite degree sequence on $n+m$ vertices such that $D_1 =(d_{1,\max}, d_{1,\max}, \ldots, d_{1,\max}, d_1,d_{1,\min}, d_{1,\min},\ldots,d_{1,\min})$ and $D_2 = (d_{2,\max}, d_{2,\max}, \ldots, d_{2,\max}, d_2,d_{2,\min}, d_{2,\min},\ldots,d_{2,\min})$, $d_{1,\min}\le d_1\le d_{1,\max}$, $d_{2,\min}\le d_2\le d_{2,\max}$. Further, let $D' = (D_1',D_2')$ be a degree sequence on $n+m$ vertices such that for all $d_i' \in D_1'$, $d_{1,\min} \le d_i'\le d_{1,\max}$, for all $d_j'\in D_2'$, $d_{2,\min}\le d_j'\le d_{2,\max}$, $\sum_{d_i \in D_1}d_i = \sum_{d_i'\in D_1'}d_i'$ and $\sum_{d_j\in D_2}d_j = \sum_{d_j'\in D_2'}d_j'$. Then $D'$ is also graphic.
        \end{sloppypar}
        \item 
        \begin{sloppypar}
        Let $D = (D_1,D_2,D_3)$ be a graphic tripartite hypergraph degree sequence on $n_1+n_2+n_3$ vertices such that for all $k=1,2,3$ $D_k =(d_{k,\max}, d_{k,\max}, \ldots, d_{k,\max}, d_k,d_{k,\min}, d_{k,\min},\ldots,d_{k,\min})$ and   $d_{k,\min}\le d_k\le d_{k,\max}$. Further, let $D' = (D_1',D_2',D_3')$ be a degree sequence on $n_1+n_2+n_3$ vertices such that for all $k=1,2,3$  and $d_k' \in D_k'$, $d_{k,\min} \le d_k'\le d_{k,\max}$ and $\sum_{d_k\in D_k}d_k = \sum_{d_k'\in D_k'}d_k'$. Then $D'$ is also graphic.
        \end{sloppypar}    
        \item 
        \begin{sloppypar}
        Let $D = (d_{\max}, d_{\max}, \ldots, d_{max}, d, d_{\min}, d_{\min},\ldots,d_{\min})$ be a graphic hypergraph degree sequence on $n$ vertices with $d_{\min}\le d\le d_{\max}$. Further let $D'$ be a hypergraph degree sequence on $n$ vertices such that for all $d'\in D'$, $d_{\min}\le d\le d_{\max}$ and $\sum_{d\in D}d = \sum_{d' \in D'}d'$. Then $D'$ is also graphic.    
        \end{sloppypar}
    \end{enumerate}
\end{theorem}
\begin{proof}
    \begin{enumerate}
      \item 
      \begin{sloppypar}
        Lemma~\ref{lem:balancing} implies that if $D = (D_1, D_2)$ is graphical, then the new sequence obtained from taking a  balancing hinge flip on $D_1$ is also graphical.
        Hence it sufficed to prove that we could obtain $D'$ from $D$ by balancing hinge flips. We first show that we could obtain $D$ from $D'$ by reverse hinge flips. WLOG, assume $D'$ is in the decreasing order:
        $$D' = (d'_{1,1}\ge \ldots\ge d'_{1,n}), (d'_{2,1}\ge \ldots \ge d_{2,m})\textrm{\ where\ } d_{1,\min}\le d'_{1,i}\le d_{1,\max}, d_{2,\min}\le d'_{2,j}\le d_{2,\max}$$
       We are going to create a sequence of degree sequences $D' = D_0, D_1, \ldots, D_m = D$ such that each $D_t = (D_{1,t},D_{2,t})$, $t= 1, 2, \ldots m$ differs from $D_{t-1} = (D_{1,t-1},D_{2,t-1})$ by a reverse hinge flip. While there are more than one degree in $D_{1,t}$ which are smaller than $d_{1,\max}$ and larger than $d_{1,\min}$, we take the largest and smallest such degrees, call them $d_i$ and $d_j$, and increase $d_i$ by $1$ and decrease $d_j$ by $1$. Any such reverse hinge flip operation either decreases the number of degrees in $D_{1,t}$  which are smaller than $d_{1,\max}$ and larger than $d_{1,\min}$ or increase the difference between the largest and smallest such degrees. Since the difference between the largest and smallest such degrees cannot be larger than $d_{1,\max}-d_{1,\min}$, in at most $d_{1,\max}-d_{1,\min}$ steps, we decrease the number of degrees which are smaller than $d_{1,\max}$ and larger than $d_{1,\min}$. Eventually, the number of such degrees will be smaller than $2$. Since we did not change the sum of the degrees, the so-obtained degree sequence will be some $D_t = (D_{1,t}, D_{2,t})$ with $D_{1,t} = D_1$. We can do the same with the other vertex class of $D'$, and eventually, we will get $D$.      
       The inverse of the inverse hinge flips will be a series of balancing hinge flips, that is, the series $D = D_m, D_{m-1}, \ldots D_0 = D'$ shows how to transform $D$ into $D'$ by balancing hinge flips.

      \end{sloppypar}
      \item 
      \begin{sloppypar}
         Lemma~\ref{lem:balancing} implies that if $D = (D_1, D_2, D_3)$ is graphical, then the new sequence obtained from taking a hinge flip on $D_1$ is also graphical.\\
         Then similar to the proof of 1, we could obtain $D'$ from $D$ by balancing hinge flips. Hence, $D$ is graphical implying that $D'$ is graphical.
      \end{sloppypar}
      \item 
      \begin{sloppypar}
         Lemma~\ref{lem:balancing} implies that if $D$ is graphical, then the new sequence obtained from taking a hinge flip on $D$ is also graphical.\\
         Then again similar to the proof of 1 and 2, we could obtain $D'$ from $D$ by balancing hinge flips. Hence, $D$ is graphical implying that $D'$ is graphical.
      \end{sloppypar}
    \end{enumerate}
\end{proof}

In general, the following well-known theorem exists on graphic bipartite degree sequences. Although the original theorem by Havel and Hakimi was given for simple graphs, it is easy to see that the following extension for bipartite graphs holds.
\begin{theorem}[Havel-Hakimi, \cite{Havel,Hakimi}]\label{theo:bipartite-HH}
Let $D = (D_1,D_2)$ be a bipartite degree sequence with $D_1 = (d_{1,1},d_{1,2},\ldots,d_{1,n_1})$ and $D_2 = (d_{2,1},d_{2,2},\ldots,d_{2,n_2})$. Then $D$ is graphic if and only if $D' = (D_1',D_2')$ with $D_1' = (d_{1,2}, d_{1,3}, \ldots,d_{1,n_1})$ and $D_2'=(d_{2,1}-1, d_{2,2}-1, \ldots, d_{2,d_{1,1}}-1, d_{2,d_{1,1}+1}, \ldots, d_{2,n_2})$ is graphic.
\end{theorem}

%\begin{theorem}[Gale-Ryser \cite{Gale,Ryser}]\label{theo:Gale-Ryser}
%Let $D = (d_{1,1} \ge d_{1,2} \ge \ldots \ge d_{1,n}), (d_{2,1},d_{2,2},\ldots, d_{2,m})$ be a bipartite degree %sequence. Then $D$ is graphic if and only if
% \begin{enumerate}
%     \item $\sum_{i=1}^n d_{1,i} = \sum_{j=1}^m d_{2,j}$ and
%     \item for all $l=1, 2, \ldots n$, $\sum_{i=1}^l d_{1,i} \le \sum_{j=1}^m \min\{l,d_{2,j}\}$.
% \end{enumerate}
%\end{theorem}

%Theorem~\ref{theo:Gale-Ryser} has a bipartite multigraph version in which we prescribe the maximal number of parallel edges.
%\begin{theorem}[\cite{Berge}]\label{theo:Berge}
%Let $D = (d_{1,1} \ge d_{1,2} \ge \ldots \ge d_{1,n}), (d_{2,1},d_{2,2},\ldots, d_{2,m})$ be a bipartite degree sequence. Then $D$ has a bipartite multigraph realization such that for all pair of vertices, $u$ and $v$, the number of parallel edges between $u$ and $v$ is at most $t$ if and only if
% \begin{enumerate}
%     \item $\sum_{i=1}^n d_{1,i} = \sum_{j=1}^m d_{2,j}$ and
%     \item for all $l=1, 2, \ldots n$, 
%     \begin{equation}
%         \sum_{i=1}^l d_{1,i} \le \sum_{j=1}^m \min\{lt,d_{2,j}\}. \label{eq:multi-gale-ryser}
%     \end{equation}
% \end{enumerate}

 By applying the Havel-Hakimi theorem, it is easy to see the following observation on almost regular bipartite degree sequences.
\begin{observation}
    Let $D= (D_1, D_2)$ be a degree sequence, such that $D_1$ is an almost regular degree sequence with some degrees $k_1$ and $k_1+1$, and $D_2$ is an almost regular degree sequence with some $k_2$ and $k_2+1$. If the degree sum of $D_1$ equals the degree sum of $D_2$, further, $k_1+1\le |D_2|$ and $k_2+1\le |D_1|$, then $D$ has a bipartite graph realization.
\end{observation}

Almost regularity will be a central concept in Section~\ref{sec:tripartite}. First, we define an almost regular configuration of incident hyperedges.
\begin{definition}
Let $H = (A, B, C, E)$ be a tripartite hypergraph, and let $v\in A$ be an arbitrary vertex. We say that the incident hyperedges of $v$ are in an almost-regular configuration if the subgraph containing only the hyperedges incident with $v$ has an almost regular degree sequence both on $B$ and $C$.
\end{definition}

It is easy to see that for any degree $d$, an almost regular configuration exists if $d \le |B|\times |C|$. Indeed, let $k_1 := \left\lfloor\frac{d}{|B|}\right\rfloor$ and $k_1' := \left\lceil\frac{d}{|B|}\right\rceil$. Obviously, $k_1\le k_1'\le |C|$ and $k_1+1 \ge k_1'$. Further, if $k_1\neq k_1'$ then there is unique $t_1$ such that 
$$
t_1\times k_1 + (|B|-t_1)\times k_1' = d.
$$
In the same way, we can define $k_2$, $k_2'$ and $t_2$ for the vertex class $C$. Then consider $D_1$ as the degree sequence containing $t_1$ $k_1$ degrees, $|B|-t_1$ $k_1'$ degrees and $D_2$ as the degree sequence containing $t_2$ $k_2$ degrees, $|C|-t_2$ $k_2'$ degrees. Since both $D_1$ and $D_2$ are both almost regular, and $d\le |B|\times |C|$, $D = (D_1, D_2)$ is a graphic bipartite degree sequence, that is, it has a realization. Consider a realization $\tilde{G} = (B,C,E)$, and let $(a,b,c)$ be a hyperedge for all $b\in B$ and $c\in C$ if and only if $(b,c) \in E(\tilde{G})$.

We also have the following observation.

\begin{observation}\label{obs:d-plus-almost-regular}
   Let $D_1 = (d, k, k, \ldots, k, k+1, k+1, \ldots, k+1)$ be a degree sequence of $n$ degrees. If $k+1\le n-1$ and $d \le |D_1|$, then $D = (D_1, D_1)$ has a bipartite graph realization.
\end{observation}
\begin{proof}
  We apply the Havel-Hakimi theorem on $D$ using the first $d$ in the first vertex class, then $d$ in the second vertex class. The remaining degree sequence is almost regular and thus, graphic.  
\end{proof}

We need a theorem for tripartite hypergraphs with almost regular degree sequences.
For this, we first prove a lemma on regular degree sequences
%%%%%%%%%%%%%%%%%%%%%%%%%%%%%%%%%%%%%%%%%
\begin{lemma}\label{lem:regular-hypergraph-lemma}
    A $k$-regular tripartite degree sequence on $n + n + n$ vertices has a tripartite hypergraph realization if $k \le n^2$.
\end{lemma}
\begin{proof}
    Assume the vertices are $(a_1,a_2,...,a_n),(b_1,b_2,...,b_n),(c_1,c_2,...,c_n)$.
    In the following construction, the indexes are in $\mathbb{Z}_n$ shifted, that is, between $1$ and $n$ modulo $n$.
    We construct the hypergraph as follows:
    For $\forall 1 \le i \le n$, we build the hyperdeges like below:\\
    $\{a_i,b_i,c_i\}, \{a_i,b_i,c_{i+1}\},...,\{a_i,b_i,c_{i+n-1}\},\\
    \{a_i,b_{i+1},c_{i+1}\},\{a_i,b_{i+1},c_{i+2}\},...,\{a_i,b_{i+1},c_{i+n}\},\\
    \ldots\\
    \{a_i,b_{i+\lfloor\frac{k}{n}\rfloor-1},c_{i+\lfloor\frac{k}{n}\rfloor-1}\},\{a_i,b_{i+\lfloor\frac{k}{n}\rfloor-1},c_{i+\lfloor\frac{k}{n}\rfloor}\},...,\{a_i,b_{i+\lfloor\frac{k}{n}\rfloor-1},c_{i+\lfloor\frac{k}{n}\rfloor+n-2}\},\\
    \{a_i,b_{i+\lfloor\frac{k}{n}\rfloor},c_{i+\lfloor\frac{k}{n}\rfloor}\}, \{a_i,b_{i+\lfloor\frac{k}{n}\rfloor},c_{i+\lfloor\frac{k}{n}\rfloor+1}\},...,\{a_i,b_{i+\lfloor\frac{k}{n}\rfloor},c_{i-\lfloor\frac{k}{n}\rfloor(n-1)+k-1}\}$\\
    
    From the construction above, we make sure that the degree of each $a_i$ is $k$ since $k\le n^2$. Indeed, if $k\le n^2$, then all the indicated hyperedges are different. Now it suffices to prove that the degree of each $b_i$ and $c_i$ is also $k$. Since each hyperedge contains one vertex in each vertex class, 
    $$\sum_{i=1}^n d(a_i) = \sum_{i=1}^n d(b_i) = \sum_{i=1}^n d(c_i).$$
    As for each $i$, $d(a_i) = k$, it suffices to prove that the degrees in the other two vertex classes are regular.
    Again from the construction, the subscripts are rotational symmetric, so $d(b_1) = d(b_2) =\ldots = d(b_n)$ and $d(c_1) = d(c_2) = \ldots = d(c_n)$.
\end{proof}
%%%%%%%%%%%%%%%%%%%%%%%%%%%%%%%%%%%%%%%%%%%%%%%%%

Now it is easy to prove the theorem on the almost regular degree sequences. 
\begin{theorem}\label{theo:hyper-almost-regular}
    Let $D_1$ be an almost regular degree sequence on $n$ vertices with $k$ or $k+1$ degrees. If 
  $k+1 \le n^2$, then $D = (D_1, D_1, D_1)$ has a tripartite hypergraph realization $H = (A,B,C,E)$.
\end{theorem}
\begin{proof}
First, we construct a tripartite hypergraph realization of the $k+1$-regular degree sequence on $n+n+n$ vertices as described in the proof of the Lemma~\ref{lem:regular-hypergraph-lemma}. For those $a_i$, $b_i$ and $c_i$ whose described degree is $k$ instead of $k+1$, we delete the hyperedge incident with these three vertices.
\end{proof}

The construction we put into the following observation will be frequently used in Section~\ref{sec:tripartite}.
\begin{observation}\label{obs:how-to-almost-regular}
    Let $D = (D_1,D_1,D_1)$ be a tripartite degree sequence with $D_1 = (d, k, k, \ldots, k, k+1, k+1, \ldots k+1)$, $|D_1|=n$ such that $d \le n^2$ and $k \ge 2n-1$. Then it is possible to exhibit hyperedges on the vertices with prescribed degree $d$ such that the remaining degrees of these vertices are $0$ and the remaining degrees of the other vertices are almost regular. By remaining degrees we mean the prescribed degrees of the vertices minus the number of incident hyperedges that we exhibit in the construction.
\end{observation}
\begin{proof}
   Let $V := \{a_0,b_0,c_0\}$ be the vertices with the prescribed degree $d$.
   In the first step, for each vertex $v \in V$, we exhibit $\min\{d,(n-1)^2\}$ hyperedges incident with $v$ and two vertices in $\bar{V}$ in an almost regular configuration with respect to $\bar{V}$. This is doable since for each vertex $v$, the hyperedges use at most $2(n-1)$ degrees, and each vertex in $\bar{V}$ has a degree at least $2n-1$. Further, since we can freely choose for which vertex in $\bar{V}$ we use $\left\lceil\frac{\min\{d,(n-1)^2\}}{n-1}\right\rceil$ and $\left\lfloor\frac{\min\{d,(n-1)^2\}}{n-1}\right\rfloor$ degrees, the remaining degrees in $\bar{V}$ can be kept almost regular.

   Next, if the remaining degree of the vertices in $V$ is $d'>0$, then for each possible combination of two vertices in $V$, we exhibit $\left\lfloor\frac{d'}{2}\right\rfloor$ hyperedges incident with the selected two vertices in $V$ and one vertex in $\bar{V}$. This is doable, since $d'\le n^2-(n-1)^2 = 2n-1$, thus $\left\lfloor\frac{d'}{2}\right\rfloor \le n-1$. Further, the remaining degrees of the vertices in $\bar{V}$ is at least $2n-1 - 2(n-1) = 1$. We can again select which vertices in $\bar{V}$ will be incident with the exhibited hyperedges, therefore the remaining degrees in $\bar{V}$ can be kept almost regular.

   Now, if there is a remaining degree of the vertices in $V$ due to parity reasons (that is, $d'$ is odd), then we add a hyperedge incident with the three vertices in $V$.

   It is clear that the hyperedges in the three steps above are disjoint.
   If we exhibit the union of the hyperedges in the three steps, the remaining degrees of the vertices in $V$ are indeed $0$, and the remaining degrees of the vertices in $\bar{V}$ are almost regular.
\end{proof}

%{\bf TODO: the following two sections must be moved to the proof of theorem 4.1.}

%Create a tripartite regular hypergraph such that the sum of the degrees is more than $\frac{2n^2}{27}$, then delete sufficiently many hyperedges from it to adjust the sum of the degrees. This causes each degree is decrease by at most cc. $\frac{\frac{2n^2}{27}}{\frac{n}{3}}$, that is $O(n)$.)\\

%From the context we know there are $\frac{n}{3}$ vertices in a vertex class. We split them into 3 vertex classes, each of which contains $\frac{n}{9}$ vertices. Then if $\frac{2n}{9} \le (\frac{n}{9})^2$,i.e.$n\ge18$, then we could apply the lemma\ref{lem:regular-hypergraph-lemma} above to construct a tripartite $\frac{2n}{9}$-regular hypergraph. Then from the context we know $O(n)-2n-4-1 = \frac{2n}{9}$ so $O(n) = \frac{20n}{9}+5$ 

%%%%%%%%%%%%%%%%%%%%%%%%%%%%%%%%%%%%%%%%%

\section{Linear bounds on graphic tripartite hypergraph degree sequences}\label{sec:tripartite}
\subsection{The $\left[\frac{2n^2}{7},\frac{5n^2}{7}\right]$ bound}
\begin{theorem}\label{theo:tripartite-bound}
 If $D = (D_1, D_2, D_3)$ is a tripartite hypergraph degree sequence on $n+n+n$ vertices such that the sums of the degrees in the three vertex classes are the same and each degree is at least $\frac{2n^2}{7}$ and at most $\frac{5n^2}{7}$, then $D$ is graphic. 
\end{theorem}
\begin{proof}
\begin{sloppypar}
    It is sufficient to prove the theorem for $D_1 = D_2 = D_3 = {\left(\5, \5, \ldots, \5, d, \2, \2, \ldots, \2\right)},$ with $\2 < d \le \5$ according to Lemma~\ref{lem:balancing}.
    For such degree sequences, we give a constructive proof. 
\end{sloppypar}
    
%   {\bf TODO: if necessary, construct the hypergraphs with a small number of vertices}

  For $n = 2$  and $3$, it is easy to verify that the theorem holds. Indeed, for $n=2$, the only possible degree sequence is $(2,2),(2,2), (2,2)$, which is graphic. Further, if $n=3$,  it is sufficient to check that $D = (D_1, D_1, D_1)$ with $D_1 = (6,6,6)$ or $(5,6,6)$  or $(4,6,6)$ $(3,6,6)$ or $(3,5,6)$ are all graphic. The remaining cases can be obtained by complementing the degree sequences (that is, subtracting each degree from $9$), constructing a realization of those degree sequences (that will be one of the above possibilities), and then complementing the obtained hypergraph. 

We call the vertices with prescribed degree $\5$ \emph{large-degree vertices}, the vertex with prescribed degree $d$ \emph{intermediate-degree vertex}, and the vertices with prescribed degree $\2$ \emph{small-degree vertices}.

    We are going to show that we can exhibit hyperedges such that the remaining degrees on the large-degree vertices (that is, the prescribed degree minus the number of incident hyperedges), as well as the remaining degrees on the small-degree vertices are almost regular. Additionally, either the remaining degrees of the intermediate-degree vertices will be $0$ (this might be done via Observation~\ref{obs:how-to-almost-regular}) or we join the intermediate-degree vertices to either the small or large-degree vertices. In the latter case, the remaining degrees of the intermediate-degree vertices will fit into the almost regular degree sequence of the chosen type of vertices.

Further, for both types of vertices, if the degrees in this remaining almost regular degree sequence are $k-1$ and $k$, then $k$ is smaller than or equal to the number of the given type vertices squared. Therefore, these almost-regular degree sequences have a tripartite hypergraph realization. Now if we take the union of the hyperedges exhibited between the small-degree and large or intermediate-degree vertices as well as the hyperedges in the realization of the almost regular degree sequences, we get a realization of $D$.

  Let $x$ be the number of intermediate and large-degree vertices in $D_1$ and let $n:=|D_1|$.
  Observe that it is sufficient to check the theorem for $x\ge \left\lceil\frac{n}{2}\right\rceil$. Indeed, if $x \le \left\lfloor\frac{n}{2}\right\rfloor$, then we can take the complement degrees, in which there will be $n-x$ large-degree vertices, thus altogether at least $n-x$ large and intermediate-degree vertices, which is at least $\left\lceil\frac{n}{2}\right\rceil$. Further, when $d < \5$, $n$ is even and $x = \frac{n}{2}$, then we can still take the complement degree sequence and get a degree sequence with $x = \frac{n}{2}+1$. That is, when $d < \5$, it is sufficient to consider $x\ge \left\lfloor\frac{n}{2}\right\rfloor+1$.  Then we can construct a realization of the complement degree sequence, take the complement of the realization and that will be a realization of the original degree sequence.
  
  We consider two cases.
\begin{itemize}[label = $\ $]
\item 
\begin{enumerate}[label= Case \arabic*:]
  \item $(x-1)^2 \ge \left\lceil\frac{2n^2}{7}\right\rceil $. First, we prove the following inequality:
   \begin{equation}
    x\times \left\lfloor\frac{5n^2}{7}\right\rfloor -  2(n-x)\times \left\lceil\frac{2n^2}{7}\right\rceil < x^3. \label{eq:x-large-inequality}
   \end{equation}
   Indeed, we divide the inequality by $n^3$ to get
   $$
    \frac{x}{n}\times\frac{ \left\lfloor\frac{5n^2}{7}\right\rfloor }{n^2}-  2\left(1-\frac{x}{n}\right)\times \frac{\left\lceil\frac{2n^2}{7}\right\rceil}{n^2} < \left(\frac{x}{n}\right)^3. 
   $$
  Now, observe that 
  $$
    \frac{ \left\lfloor\frac{5n^2}{7}\right\rfloor }{n^2} \le \frac{5}{7}
  $$ 
  and
  $$
     \frac{\left\lceil\frac{2n^2}{7}\right\rceil}{n^2} \ge \frac{2}{7},
  $$
  thus
  $$ 
      \frac{x}{n}\times\frac{ \left\lfloor\frac{5n^2}{7}\right\rfloor }{n^2}-  2\left(1-\frac{x}{n}\right)\times \frac{\left\lceil\frac{2n^2}{7}\right\rceil}{n^2} \le   \frac{x}{n}\times \frac{5}{7} - 2 \left(1-\frac{x}{n}\right)\times \frac{2}{7}.
  $$
  Now introducing the new variable $z:= \frac{x}{n}$, we have the inequality
 $$
  \frac{5}{7}z-2(1-z)\frac{2}{7} < z^3
 $$
  which holds for all positive $z$, thus the inequality in equation~\ref{eq:x-large-inequality} also holds.

  The inequality in equation~\ref{eq:x-large-inequality} tells the following: if we are able to exhibit $\left\lceil\frac{2n^2}{7}\right\rceil$ hyperedges on each small-degree vertex and the other two incident vertices of each such hyperedge are not small-degree vertices, then the remaining degree sum of the large and intermediate-degree vertices is less than $x^3$.

  We construct the appropriate hyperedges in the following way.
  We index thesmall-degree vertices in each vertex class, and for each $i= 1, 2, \ldots n-x$, we consider the set $V_i = \{a_i,b_i,c_i\}$.
  For each $i = 1, 2, \ldots n-x$ and for each vertex $v\in V_i$, we exhibit $\left\lceil\frac{2n^2}{7}\right\rceil$ hyperedges incident with one $v$ and two intermediate or large-degree vertices. When we talk about the remaining degrees of the large or intermediate-degree vertices after exhibiting hyperedges incident with vertices from $V_i$, we mean their prescribed degrees minus the number of the incident hyperedges exhibited altogether for all $V_j$, $j=1,2,\ldots, i$.

  We exhibit the hyperedges in the following way. While the remaining degrees of the large-degree vertices are larger than $d$, we exhibit hyperedges in an almost regular configuration on the large-degree vertices, such that the remaining degrees of the large-degree vertices are almost regular. This is doable according to Observation~\ref{obs:how-to-almost-regular} since $(x-1)^2 \ge \left\lceil\frac{2n^2}{7}\right\rceil $. There might exist a $V_i$ such that after exhibiting hyperedges incident with vertices in it, the remaining degrees of the large-degree vertices become smaller than $d$. Then instead of exhibiting the almost regular configuration on the large-degree vertices, we prescribe degrees on the intermediate and large-degree vertices such that the remaining degree sequence becomes almost regular. Let the remaining degrees of the large-degree vertices before exhibiting hyperedges incident with vertices in $V_i$ be almost regular for some $k$'s and $(k+1)$'s, and we define 
  $f:= \2-(k+1-d)(x-1)$. Observe that $f > 0$ since the remaining degrees of the large-degree vertices would be smaller than $d$ if the hyperedges were exhibited in an almost regular configuration on the large-degree vertices. 
  Then we are seeking a configuration (a bipartite graph realization) of the degree sequence $\tilde{D} =(\tilde{D}_1,\tilde{D}_1)$ with
  $$
\tilde{D}_1:=  \left(\left\lceil\frac{f}{x}\right\rceil, k+1-d+\left\lceil\frac{f}{x}\right\rceil, \ldots,k+1-d+\left\lceil\frac{f}{x}\right\rceil, k+1-d+\left\lfloor\frac{f}{x}\right\rfloor, \ldots k+1-d+\left\lfloor\frac{f}{x}\right\rfloor\right)$$
  with the appropriate numbers of ceiled and floored values, such that the degree sum is $\2$.
  Such a configuration exists according to Observation~\ref{obs:d-plus-almost-regular}.
  When the remaining degrees of the large and intermediate-degree vertices start forming an almost regular degree sequence, we exhibit hyperedges in an almost regular configuration on them according to Observation~\ref{obs:how-to-almost-regular}. We claim this procedure will not terminate before adding hyperedges to all small-degree vertices. This is because 
  $$
  2\times \left\lceil\frac{2n^2}{7}\right\rceil \times (n-x) \le d + \left\lfloor\frac{5n^2}{7}\right\rfloor\times (x-1).
  $$
  Indeed, it is easy to see that whenever $(x-1)^2 \ge \left\lceil\frac{2n^2}{7}\right\rceil$, $x\ge \frac{n+1}{2}$. Therefore, $n-x \le x-1$. Further, if $n>2$, then $ 2\times \left\lceil\frac{2n^2}{7}\right\rceil \le \left\lfloor\frac{5n^2}{7}\right\rfloor$. %(For $n=2$, $x$ should be larger than $2$ to satisfy the inequality of Case 1, an impossible setup).
 
  Since we do this procedure for each small-degree vertex in all three vertex classes, the remaining degree sum of the large and intermediate-degree vertices in the three vertex classes will be the same. Further, the remaining degrees of the intermediate-degree vertices will be the same in the three vertex classes. In addition, in all three vertex classes, the remaining degrees of the large-degree vertices are almost regular, and the degree sequences are the same (the order of them might be different, but this is not an issue).

  Now we claim that the remaining degrees of the intermediate-degree vertices are smaller than $x^2$. Indeed, if it was larger than $x^2$, then all degrees of the $x$ intermediate and large-degree vertices in $D_1$  (as well as in $D_2$ and $D_3$) were larger than  $x^2$, then the remaining degrees altogether were more than $x^3$, contradicting the proved inequality in equation~\ref{eq:x-large-inequality}.

It is clear that the remaining degrees of the small-degree vertices are $0$ since we exhibited $\2$ hyperedges incident with each of these vertices. There are two possibilities for the remaining degrees of the large and intermediate-degree vertices. Either we reached the point where we started exhibiting hyperedges incident with the intermediate-degree vertices or not. If we reach the point, then the remaining degrees of the intermediate and large-degree vertices are almost regular, and the maximum degree is at most $x^2$ (according to the inequality in equation~\ref{eq:x-large-inequality}), therefore, they have a tripartite hypergraph realization. The union of this tripartite hypergraph realization and the exhibited hyperedges incident with one small-degree vertex form a realization of the prescribed degree sequence.

If we did not reach the point where we exhibit hyperedges incident with intermediate-degree vertices, then the remaining degree of the intermediate-degree vertex is $d$. It is clear that $d<x^2$ since the remaining degrees of the large-degree vertices are larger than $d$, and this would contradict the inequality in equation~\ref{eq:x-large-inequality}. Since $d\le x^2$, we can introduce the following hyperedges incident with at least one intermediate-degree vertex and no small-degree vertex. First, for each intermediate-vertex, we add $\min\{d, (x-1)^2\}$ hyperedges incident with one intermediate and two large-degree vertices, in an almost regular configuration. Then if the remaining degree is $d'>0$, for all possible combinations of two intermediate and one large-degree vertex with respect to the three vertex classes, we add $\left\lfloor\frac{d'}{2}\right\rfloor$ hyperedges in an almost regular configuration. If $d'$ is odd and there is one remaining degree on each of the intermediate-degree vertices, we add the hyperedge incident with the three intermediate-degree vertices. This is all doable as $d<x^2$ and the remaining degrees of the large-degree vertices are larger than $d$. Now we claim that the remaining degrees of the large-degree vertices are at most $(x-1)^2$. Indeed, we can treat the intermediate-degree vertices as the small-degree vertices and replace $x$ with $x' = x-1$ in the inequality~\ref{eq:x-large-inequality}. The remaining degrees of the large-degree vertices are at most $(x-1)^2$ and they form an almost regular degree sequence, therefore, there is a tripartite hypergraph realization of them. The union of the tripartite hypergraph realization on the intermediate and large-degree vertices and the exhibited hyperedges incident with one small-degree vertex forms a realization of the prescribed degree sequence.

  \item $\left\lceil\frac{n}{2}\right\rceil \le x$ and $(x-1)^2 < \left\lceil\frac{2n^2}{7}\right\rceil$. The challenge here is how to spread the hyperedges incident with one small-degree and two large or intermediate-degree vertices as well as how to treat the intermediate-degree vertices. First, we discuss the subcase when $d = \5$, and then the subcase when $d< \5$.

\noindent{\it Subcase 2a: $d=\5$}. In this case, there is no intermediate-degree vertex. Observe that the maximum number of hyperedges incident with one small-degree vertex and two large-degree vertices incident with one particular small-degree vertex is $\min\{\2,x^2\}$.  We again index the small-degree vertices in each vertex class, and for each $i= 1, 2, \ldots n-x$, we consider the set $V_i = \{a_i,b_i,c_i\}$.
  For each $i = 1, 2, \ldots n-x$ and for each vertex $v\in V_i$, we exhibit $\min\{\2,x^2\}$ hyperedges incident with one $v$ and two large-degree vertices in an almost regular configuration with respect to the large-degree vertices.

We need to show three things.
\begin{enumerate}
\item The remaining degrees of the large-degree vertices do not drop below $0$. Indeed, each hyperedge incident with one small-degree vertex and two large-degree vertices takes one degree from the small-degree vertex and two degrees from large-degree vertices. Since the maximum number of hyperedges in question is $(n-x)\times \2$, and the degree sum of the large-degree vertices is $x\times \5$, we need to show that
$$
2\times (n-x)\times \2 \le x\times \5.
$$
Indeed, $n-x < x$ and $2\times \2 < \5$ for any $n>3$, therefore the inequality holds.
\item The remaining degrees of the large-degree vertices are almost regular and are smaller than or equal to $x^2$. 
Indeed, we exhibit $\min\left\{2x^2(n-x),2\left\lceil\frac{2n^2}{7}\right\rceil(n-x) \right\}$ hyperedges on the large-degree vertices altogether, depending if $x^2$ or $\2$ hyperedges are exhibited on the small-degree vertices. If the minimum is taken at $2x^2(n-x)$, we have to prove that
  \begin{equation}
  \left\lfloor\frac{5n^2}{7}\right\rfloor  x - 2x^2(n-x) \le x^3.\label{eq:large-becomes-small-enough}
  \end{equation}
  Observe that $\frac{\left\lfloor\frac{5n^2}{7}\right\rfloor}{n^2} \le \frac{5}{7}$. Dividing by $n^2x$ and introducing the new variable $z := \frac{x}{n}$, we get that
  $$
  \frac{5}{7} - 2z(1-z) \le z^2
  $$
  which holds for any $z$ between $0.5$ and $1$. 

  If the minimum is taken at $2\left\lceil\frac{2n^2}{7}\right\rceil(n-x)$, then we have to prove that
  $$
  \left\lfloor\frac{5n^2}{7}\right\rfloor x- 2\left\lceil\frac{2n^2}{7}\right\rceil(n-x) \le x^3.
  $$
  Observing that $\frac{\left\lceil\frac{2n^2}{7}\right\rceil}{n^2} \ge \frac{2}{7}$, dividing by $n^3$ and introducing the new variable $z := \frac{x}{n}$, we have to prove that
  \begin{equation}
  \frac{5}{7}z-\frac{4}{7}(1-z) \le z^3, \label{eq:limit}      
  \end{equation}
  
  that holds for any positive $z$.

\item The remaining degrees of the small-degree vertices are at most $(n-x)^2$. We have to prove that whenever $\left\lceil\frac{2n^2}{7}\right\rceil > x^2$,
  \begin{equation}
  \left\lceil\frac{2n^2}{7}\right\rceil- x^2 \le (n-x)^2.\label{eq:small-is-small-enough}
  \end{equation}
  Observe that $\left\lceil\frac{2n^2}{7}\right\rceil \le \frac{2n^2}{7}+1$. Then dividing the inequality by $n^2$, introducing the new variable $z:= \frac{x}{n}$ and assuming that $n\ge 3$, we get that
  $$
  \frac{\left\lceil\frac{2n^2}{7}\right\rceil}{n^2}- z^2 \le \frac{2}{7}+\frac{1}{n^2}-z^2 \le \frac{2}{7}+\frac{1}{9}-z^2 \le (1-z)^2.
  $$
  The last inequality holds for any $z$. 
\end{enumerate}
Therefore, the remaining degrees of the large-degree vertices are almost regular and at most $x^2$. They have a tripartite hypergraph realization. Similarly, the remaining degrees of the small-degree vertices are regular and at most $(n-x)^2$. Therefore, they also have a tripartite hypergraph realization. These three types of hyperedges (hyperedges incident with one small-degree vertex and two large-degree vertices, hyperedges incident with three small-degree vertices, and hyperedges incident with three large-degree vertices) together form a realization of the prescribed degree sequence.

  \vspace{0.3cm}
  \noindent{\it Subcase 2b: $d<\5$.} This subcase is split further into two sub-subcases.

\noindent{\it Sub-subcase 2b i):  $d\le x^2$} In this sub-subcase, we handle the intermediate-degree vertex as if it were a small-degree vertex. That is, we first exhibit for each intermediate-degree vertex $d^*:= d-(\2-(x-1)^2)$ hyperedges incident with intermediate-degree vertices and at most two large-degree vertices. Observe that $d^*$ is positive, and actually, larger than $(x-1)^2$ as it is $d- \2 + (x-1)^2$ and $d>\2$. Also observe that it is less than $x^2$, since $\2 > (x-1)^2$. 

We do it in two phases. First, for each intermediate-degree vertex $v$, we exhibit $(x-1)^2$  hyperedges incident with $v$ and two large-degree vertices in a(n almost) regular configuration with respect to the large-degree vertices. Then let the number of hyperedges we want to add further to be incident with each of the intermediate-degree vertices be $d':= d^* - (x-1)^2$. We exhibit $\left\lfloor\frac{d'}{2}\right\rfloor$ hyperedges in each possible combination with respect to the three vertex classes incident with two intermediate and one large-degree vertices. Observe that there is at most one hyperedge of these types incident with a particular large-degree vertex. Indeed, $d^* < x^2$ for $d\le x^2$. Therefore, $\left\lfloor\frac{d'}{2}\right\rfloor = \left\lfloor\frac{1}{2}*(d^* - (x-1)^2)\right\rfloor < x-\frac{1}{2}$ and there are $x$ large-degree vertices in each degree class.
If $d'$ is odd, then $d'-2\left\lfloor\frac{d'}{2}\right\rfloor = 1$.
In this case, we exhibit a hyperedge incident with three intermediate-degree vertices. Then the remaining degree of the intermediate-degree vertices is $\2-(x-1)^2$ and the remaining degrees of the large-degree vertices are almost regular.

Next for each small-degree vertex, we exhibit $(x-1)^2$  hyperedges incident with one small-degree vertex and two large-degree vertices. Now the remaining degrees of the small-degree vertices and the intermediate-degree vertices are $\2-(x-1)^2$-regular on $n-x+1$ vertices. If $n$ is odd and $x = \left\lceil\frac{n}{2}\right\rceil$, then we add further hyperedges incident with two small or intermediate-degree and one large-degree vertices. We need to do this because if $n$ is small, then the remaining degrees of the large-degree vertices could be larger than $(x-1)^2$. For each large-degree vertex with remaining degree $d"$, we exhibit $\min\left\{d",\left\lfloor\frac{\2-(x-1)^2}{2}\right\rfloor\right\}$ hyperedges in an almost regular configuration with respect to the small and intermediate-degree vertices. Observe that it cannot bring the remaining degrees of the small and intermediate-degree vertices below $0$. Indeed, before adding hyperedges incident with one large and two small or intermediate-degree vertices, the remaining degrees of the small and intermediate-degree vertices were $\2-(x-1)^2$. The maximum number of hyperedges incident with one large-degree and two small or intermediate-degree vertices is 
$$
3\times (x-1) \times \left\lfloor\frac{\2-(x-1)^2}{2}\right\rfloor.
$$
Each hyperedge takes two degrees from the small or intermediate-degree vertices. However, the degree sum of these vertices is $3\times(n-x+1)\times\left(\2-(x-1)^2)\right)$ and
$$
2\times 3\times (x-1) \times \left\lfloor\frac{\2-(x-1)^2}{2}\right\rfloor < 3\times(n-x+1)\times\left(\2-(x-1)^2\right),
$$
since $x =  \left\lceil\frac{n}{2}\right\rceil$ and thus $x-1 < n-x+1$.

The remaining degrees of the intermediate and small-degree vertices will be at most $\2-(x-1)^2$ and almost regular on $n-x+1$ vertices. Therefore, it will be a graphic degree sequence. This can be proved by replacing $x$ with $x-1$ in equation~\ref{eq:small-is-small-enough}.
Yet we need to show two things.
\begin{enumerate}
    \item The remaining degrees of the large-degree vertices remain non-negative. For this, we need to show that 
$$
2\times (x-1)^2\times (n-x+1) + (x-1) \le (x-1) \times \5.
$$
Dividing by $(x-1)$ we get
$$
2\times (x-1)\times (n-x+1) + 1 \le \5,
$$
which always holds since the left-hand side is smaller than or equal to $\frac{n^2}{2}+1$.
  \item The remaining degrees of the large-degree vertices are smaller than $(x-1)^2$. If $x-1$ is still at least $\left\lceil\frac{n}{2}\right\rceil$, then this can be proven by replacing $x$ with $x-1$ in equation~\ref{eq:large-becomes-small-enough}. $x-1$ can become smaller than $\left\lceil\frac{n}{2}\right\rceil$ only if $n$ is odd and $x=\left\lceil\frac{n}{2}\right\rceil$. We carefully handle the case by exhibiting further $\min\left\{d",\left\lfloor\frac{\2-(x-1)^2}{2}\right\rfloor\right\}$ hyperedges incident with the large-degree vertices. This makes the remaining degrees of the large-degree vertices either $0$ or at most
$$
\5-2\left\lfloor\frac{n}{2}\right\rfloor^2- \left\lfloor\frac{\2-\left\lfloor\frac{n}{2}\right\rfloor^2}{2}\right\rfloor
$$ 
It converges to $\frac{11}{56}n^2$ which is smaller than $(x-1)^2 \approx \frac{n^2}{4}$. For small odd $n$-s, it is easy to verify that the remaining degree of the large-degree vertices is always smaller than $(x-1)^2$.
\end{enumerate}

\noindent{\it Sub-subcase 2b ii):  $d> x^2$} In this case, we treat the intermediate-degree vertices as if they were large-degree vertices. However, when adding hyperedges incident with one small-degree and two intermediate or large-degree vertices, we design configurations such that the remaining degrees of the intermediate-degree vertices will be exactly $x^2$ if $d-x^2$ is even, and $x^2-1$ if $d-x^2$ is odd. That is, let $p := (d-x^2)\mod 2$, $d^* := d-x^2+p$, and let $\bar{f} := \left\lceil\frac{d^*}{2(n-x)}\right\rceil$, further $\underline{f} := \left\lfloor\frac{d^*}{2(n-x)}\right\rfloor$. Let $t$ be the unique integer for which
$$
t\times \underline{f} + (n-x-t)\times \bar{f} = \left\lfloor\frac{d^*}{2}\right\rfloor = \frac{d^*}{2}.
$$
Observe that $d^*$ is always an even number, so the last equality holds.

We exhibit $\min\{\2,\underline{f}+(x-1)^2\}$ hyperedges on each of $t$ small-degree vertices which are incident with one small-degree vertex and two intermediate or large-degree vertices, such that in the configuration with respect to the intermediate and large-degree vertices, there is a degree $\underline{f}$ on the intermediate-vertex and the degrees are almost regular on the large-degree vertices. Also, we exhibit $\min\{\2,\bar{f}+(x-1)^2\}$ hyperedges on each of $n-x-t$ small-degree vertices which are incident with one small-degree vertex and two intermediate or large-degree vertices, such that in the configuration with respect to the intermediate and large-degree vertices, there is a degree $\bar{f}$ on the intermediate-vertex and the degrees are almost regular on the large-degree vertices. First, we need to show that $\bar{f}$ is at most $x$. Indeed, assume for contradiction that $\bar{f} > x$. Because $\bar{f}$ and $x$ are integers, we have $\bar{f}-1 \ge x$. In addition, it is easy to see $\underline{f} \ge \bar{f}-1$ so $\underline{f} \ge x$. Then
$$
\frac{d - x^2+1}{2(n-x)} \ge\left\lfloor\frac{d - x^2+1}{2(n-x)}\right\rfloor \ge \left\lfloor\frac{d^*}{2(n-x)}\right\rfloor \ge x
$$
and then it would imply that 
$$
\frac{5n^2}{7} \ge \5 \ge d+1\ge x^2 + 2(n-x)x
$$
Dividing it by $n^2$ and introducing the new variable $z:= \frac{x}{n}$, we would get
$$
\frac{5}{7} \ge z^2 + 2(1-z)z,
$$
however, it does not hold for any $0.5 \le z \le 1$. Therefore, $\bar{f}\le x$, and thus the requested configuration exists.

If $\min\{\2,\bar{f}+(x-1)^2\} = \2$, then the degree sequence of the configuration takes $\bar{f}$ degree on the intermediate-degree vertex, and an almost regular configuration on the large-degree vertices. It is possible since $\bar{f} \le x \le \2$. Same holds for $\min\{\2,\underline{f}+(x-1)^2\} = \2$.

After adding the hyperedges incident with one small-degree vertex, we add the remaining $x^2$ or $x^2-1$ hyperedges incident with the intermediate-degree vertices and at most two large-degree vertices. First, we add $(x-1)^2$ hyperedges incident with one intermediate-degree vertex and two large-degree vertices in an (almost) regular configuration. Then for each combination of two intermediate-degree and one large-degree vertices with respect to the vertex classes, we add $x-1$ hyperedges incident with two intermediate and one large-degree vertices. If the prescribed remaining degree of the intermediate-degree vertices was $x^2$, then we also add the hyperedge incident with the three intermediate-degree vertices. Now the remaining degree of the intermediate-degree vertices is $0$.

We also need to show three things.
\begin{enumerate}
    \item The remaining degrees of the large-degree vertices are non-negative. We subtracted at most $2(n-x)x$ degrees by hyperedges incident with one small-degree vertex and two intermediate or large-degree vertices, and at most $2x-1$ degrees by hyperedges incident with at least one large-degree and at least one and at most two intermediate-degree vertices. That is, we need that
    $$
    2(n-x)x +2x-1  = -2x^2+2(n+1)x -1\le \5.
    $$
    However, the left hand side is at most $\frac{n^2}{2}+n-\frac{1}{2}$, and thus the inequality always holds for all $n>5$. For $n=4$ and $n=5$, a manual investigation of all possible cases reveals that the remaining degrees of the large-degree vertices are always non-negative.
    \item The remaining degrees of the large-degree vertices are at most $(x-1)^2$. There are two cases: $\min\{\2,\bar{f}+(x-1)^2\} = \2$ or $\min\{\2,\bar{f}+(x-1)^2\} = \bar{f}+(x-1)^2$. If $\min\{\2,\bar{f}+(x-1)^2\} = \2$, then observe that the subtracted degrees on the large-degree vertices in sub-subcase 2b ii) cannot be less than those in Subcase 2a. Indeed, if $d$ was $\5$, then the remaining degrees of the large-degree vertices were almost regular, and as we proved, each degree is at most $x^2$. Now, we also subtract as many degrees as possible ($\2$ with each small-degree vertex), however, the remaining degrees of the intermediate-degree vertices are $x^2$ or $x^2-1$, thus the other degrees cannot be larger. Further, then we took $2x-1$ degrees from each large-degree vertex by exhibiting hyperedges on one or two intermediate-degree and one or two large-degree vertices. 

On the other hand, if $\min\{\2,\bar{f}+(x-1)^2\} = \bar{f}+(x-1)^2$, then the remaining degrees of the large-degree vertices cannot be larger than those in sub-subcase 2b i). Indeed, in sub-subcase 2b i), we subtract $x-1$ degree with each small-degree vertex, and here we subtract $x-1$ too. In sub-subcase 2b i), we subtract at most $2x-1$ degree from the large-degree vertices using the intermediate-degree vertices. However, here we subtract exactly $2x-1$. Therefore, the remaining degrees of the large-degree vertices in sub-subcase 2b ii) are also at most $(x-1)^2$.
    \item We also need to show that the remaining degrees of the small-degree vertices are at most $(n-x)^2$. We subtract at least $(x-1)^2$ degrees from them. Therefore, we need to show that
$$
\2 - (x-1)^2 \le (n-x)^2
$$
  Observe that $\left\lceil\frac{2n^2}{7}\right\rceil \le \frac{2n^2}{7}+1$. Then dividing the inequality by $n^2$, introducing the new variable $z:= \frac{x}{n}$, we get that
  $$
  \frac{\left\lceil\frac{2n^2}{7}\right\rceil}{n^2}- \left(z-\frac{1}{n}\right)^2 \le \frac{2}{7}+\frac{1}{n^2}-\left(z-\frac{1}{n}\right)^2 \le (1-z)^2.
  $$
  The last inequality holds for any $z$ if $n\ge 6$. For $n=4$ and $n=5$, it can be manually checked that the remaining degrees of the small-degree vertices are always smaller than or equal to $(n-x)^2$.

\end{enumerate}

  Therefore, on the remaining degrees, there is a tripartite hypergraph realization on both the small-degree and the large-degree vertices. If we take the union of all the exhibited hyperedges, we get a realization of the prescribed degree sequence.

 \end{enumerate}
\end{itemize}
\end{proof}

\subsection{Algorithmic considerations}\label{ssec:algorithmic-tripartite}

In this subsection, we explain how the constructive proof of Theorem~\ref{theo:tripartite-bound} provides a polynomial running time algorithm. We do not give a thorough analysis but just the reasons why we can construct a tripartite hypergraph realization of a degree sequence on $n+n+n$ vertices bounded by $\frac{2n^2}{7}$ and $\frac{5n^2}{7}$ in polynomial time.

Let $(D_1, D_2, D_3)$ be a hypergraph degree sequence on $n+n+n$ points satisfying the conditions of Theorem~\ref{theo:tripartite-bound}. Let $\Sigma$ denote the degree sum in each vertex class. It is easy to find the integers $k$ and $d$ such that
$$
k\times \frac{5n^2}{7} + (n-k-1)\frac{2n^2}{7} + d = \Sigma.
$$
Let $\tilde{D}$ be the degree sequence that contains $k$ $\frac{5n^2}{7}$'s, $(n-k-1)$ $\frac{2n^2}{7}$'s, and one $d$ ($d$ could be $\frac{5n^2}{7}$ or $\frac{2n^2}{7}$). First, we construct a tripartite hypergraph realization of $(\tilde{D},\tilde{D},\tilde{D})$, and then we apply the procedure in the proof of  Theorem~\ref{theo:balancing} to get a realization of $(D_1, D_2, D_3)$. 

To construct a hypergraph realization of $(\tilde{D},\tilde{D},\tilde{D})$, we just need to follow the explicit construction in the proof of Theorem~\ref{theo:tripartite-bound}, which clearly can be done in polynomial time. 
The proof of Theorem~\ref{theo:balancing} gives an algorithm on how to construct a series of degree sequences from $(\tilde{D},\tilde{D},\tilde{D})$ to $(D_1, D_2, D_3)$ such that each consecutive pair of degree sequences differ by a balancing hinge flip operation.
Each balancing hinge flip operation in Theorem~\ref{theo:balancing} decreases the $L_1$ distance between the current degree sequence and the prescribed degree sequence. Since the difference is an integer, the degree sum is $O(n^3)$. A balancing hinge flip operation on the degree sequence and its corresponding realization can be found in polynomial time. Constructing a realization of $(D_1, D_2, D_3)$ from a realization of $(\tilde{D},\tilde{D},\tilde{D})$ clearly can be done in polynomial time. Therefore, constructing a realization of $(D_1, D_2, D_3)$ can be done in polynomial time.

\subsection{Discussion on the limits}
It is meaningful to ask what is the smallest $c$ such that the tripartite hypergraph degree sequences on $n+n+n$ vertices with degrees between $cn^2$ and $(1-c)n^2$ are always graphic. We feel that equation~\ref{eq:limit} plays a central role. Let $c$ be the real number between $0$ and $1$ such that the equation
\begin{equation}
    (1-c)z - 2c(1-z) = z^3\label{eq:tripartite-bound-c}
\end{equation}
has a double root. The only positive real solution for $c$ is
\begin{equation}
%    8 - 3/2 (1 + i sqrt(3)) (1/2 (37 + i sqrt(3)))^(1/3) - (21 (1 - i \sqrt(3)))/(2^(2/3) (37 + i \sqrt(3))^(1/3)
c = 8 - \frac{3}{2} (1 + i \sqrt{3})\sqrt[3]{\frac{1}{2} (37 + i \sqrt{3})} - \frac{21 (1 - i \sqrt{3})}{\sqrt[3]{4 (37 + i \sqrt(3)}}\approx 0.278066.\label{eq:c}
\end{equation}
We have the following conjecture.
\begin{conjecture}\label{con:tripartite-bound}
  For any $\varepsilon > 0$, there exists an $n_0$ such that for any $n\ge n_0$, any tripartite hypergraph degree sequence on $n+n+n$ vertices with equal degree sums in the three vertex classes, and all degrees between $(c+\varepsilon)n^2$ and $(1-c-\varepsilon)n^2$ with the $c$ in equation~\ref{eq:c} is graphic.

  Further, for any $\varepsilon > 0$, there exists an $n_0$ such that for any $n\ge n_0$, there exists a tripartite hypergraph degree sequence on $n+n+n$ vertices such that the degree sums in the three vertex classes are equal, all degrees are between $(c-\varepsilon)n^2$ and $(1-c+\varepsilon)n^2$, but the tripartite hypergraph degree sequence is not graphic.
\end{conjecture}

Observe that $\frac{2}{7}\approx 0.2857$, so Theorem~\ref{theo:tripartite-bound} is close to the conjecture above. Notice that proving Theorem~\ref{theo:tripartite-bound} already includes many technical difficulties, and these technical difficulties might increase rapidly as we get closer to the conjectured bound $c$. Also, due to some rounding issues, some tripartite hypergraph degree sequences within the prescribed bounds but on a small number of vertices might become non-graphic.

It is easy to prove the following observation using the Gale-Ryser theorem \cite{Gale,Ryser} and Theorem~\ref{theo:balancing}.
\begin{observation}\label{obs:bipartite-bound}
   Any bipartite degree sequence on $n+n$ vertices, with the same degree sums on the two vertex classes, and each degree between $\frac{n}{4}$ and $\frac{3n}{4}$ is graphic. 
\end{observation}
For tripartite hypergraphs, the conjectured bounds are slightly narrower. Still, it is easy to prove the following observation.
\begin{observation}
    Let $n$ be even, then any tripartite hypergraph degree sequence on $n+n+n$ vertices with the degree sum $\frac{n^3}{2}$ in each vertex class, and all degrees between $\frac{n^2}{4}$ and $\frac{3n^2}{4}$ is graphic.
\end{observation}
\begin{proof}
    According to Theorem~\ref{theo:balancing}, it is sufficient to check that $D = (D_1,D_1,D_1)$ with $D_1$ containing $\frac{n}{2}$ $\frac{n^2}{4}$'s and $\frac{n}{2}$ $\frac{3n^2}{4}$'s is graphic. It is easy to construct a realization of $D$. Indeed, add all hyperedges incident with one vertex with the prescribed degree $\frac{n^2}{4}$ and two vertices with the prescribed degree $\frac{3n^2}{4}$, and further add all hyperedges incident with three vertices with the prescribed degree $\frac{3n^2}{4}$. It is easy to verify that this is indeed a realization of $D$.
\end{proof}
On the other hand, it is easy to see that $D = (D_1, D_1, D_1)$ with $D_1 = (9,9,27,27,27,27)$ is not graphic ($9= \frac{6^2}{4}$, $27=\frac{3*6^2}{4}$). However, $D =(D_1,D_1,D_1)$ with $D_1 = (9,27,27,27,27,27)$ is again graphic. If we consider equation~\ref{eq:tripartite-bound-c} with $c= \frac{1}{4}$, we get that there are two positive solutions for $z$:
$$
z_1 = 0.5,\ \ \ \ \ z_2 = \frac{\sqrt{17}-1}{4} \approx 0.78078. 
$$
Based on these, we have the following conjecture.
\begin{conjecture}
 For any $\varepsilon>0$, there exists an $n_0$ such that for all $n>n_0$, any tripartite hypergraph degree sequence on $n+n+n$ vertices, with each degree between $\frac{n^2}{4}$ and $\frac{3n^2}{4}$, the degree sums in the three vertex classes equal, and the degree sum is either at most $\left(\frac{5-\sqrt{17}}{4}-\varepsilon\right)\frac{3n^3}{4}+ \left(\frac{\sqrt{17}-1}{4}+\varepsilon\right)\frac{n^3}{4} = \left(\frac{7-\sqrt{17}-\varepsilon}{8}\right)n^3\left[\approx (0.36-\varepsilon)n^3\right]$ or at least $\left(\frac{\sqrt{17}-1}{4}+\varepsilon\right)\frac{3n^3}{4}+ \left(\frac{5-\sqrt{17}}{4}-\varepsilon\right)\frac{n^3}{4}=\left(\frac{1+\sqrt{17}+\varepsilon}{8}\right)n^3\left[\approx (0.64+\varepsilon)n^3\right]$ simultaneously in each vertex class is graphic.

 On the other hand, for any $\frac{\sqrt{17}-3}{4}>\varepsilon>0$, there exists an $n_0$ such that for all $n>n_0$, there exists a tripartite hypergraph degree sequence $D$ on $n+n+n$ vertices, such that the degree sum in each vertex class is $\left\lceil\left(\frac{5-\sqrt{17}}{4}+\varepsilon\right)\frac{3n^3}{4}+ \left(\frac{\sqrt{17}-1}{4}-\varepsilon\right)\frac{n^3}{4}\right\rceil$, each degree is between $\frac{n^2}{4}$ and $\frac{3n^2}{4}$ and $D$ is not graphic.
\end{conjecture}
[Remark: we need the upper bound for $\varepsilon$ in the conjecture to avoid the case when the degree sum is exactly $\frac{n^3}{2}$ or at least $\left(\frac{\sqrt{17}-1}{4}+\varepsilon\right)\frac{3n^3}{4}+ \left(\frac{5-\sqrt{17}}{4}-\varepsilon\right)\frac{n^3}{4}$.]

\section{Linear bounds on hypergraph degree sequences}\label{sec:linear-hypergraphic}
\subsection{The $[\frac{2n^2}{63}+\frac{4n}{5}+5, \frac{5n^2}{63} - \frac{20n}{63}]$ bound}
In this section, we prove that any hypergraph degree sequence on $n$ vertices is graphic if the sum of the degrees can be divided by $3$ and all degrees are between $\frac{2n^2}{63}+O(n)$ and $\frac{5n^2}{63}-O(n)$.
The idea is to construct an ``almost" tripartite hypergraph and a few additional hyperedges that span over two vertex classes or inside one vertex class. The additional hyperedges use $O(n)$ degrees from each vertex, such that we arrive at a tripartite hypergraph degree sequence for which we can apply Theorem~\ref{theo:tripartite-bound}. %{\bf (TODO: or some other bounds)}
\begin{theorem}\label{theo:hypergraph-linear-bound}
If $D$ is a hypergraph degree sequence on $n$ vertices such that each degree is at least  $\frac{2\left\lfloor\frac{n}{3}\right\rfloor^2}{7}+4*\left\lceil\frac{n}{5}\right\rceil+1$ and at most $\frac{5\left\lfloor\frac{n}{3}\right\rfloor^2}{7}$, $n\ge 45$ and the sum of the degrees can be divided by $3$, then $D$ is graphic.
\end{theorem}
Before the proof, we add a remark. It is easy to see the following. If $d \ge \frac{2n^2}{63}+\frac{4n}{5}+5$, then $d \ge \frac{2\left\lfloor\frac{n}{3}\right\rfloor^2}{7}+4*\left\lceil\frac{n}{5}\right\rceil+1$ whenever $n>1$. Further, if $d \le \frac{5n^2}{63} - \frac{20n}{63}$ then $d \le \frac{5\left\lfloor\frac{n}{3}\right\rfloor^2}{7}$. Therefore, we could write a more readable, though slightly narrower bound that the degrees should be between $\frac{2n^2}{63}+\frac{4n}{9}+5$ and $\frac{5n^2}{63} - \frac{20n}{63}$.
\begin{proof}

We build up a realization using four types of hyperedges. We add the hyperedges by types, that is, we build up a realization in four phases.  After adding the hyperedges of a given type (that is, after each phase), we adjust the degrees of the vertices by subtracting the number of incident hyperedges of the given type from their prescribed degrees. Then in the next phase, we consider these degree sequences. In the last phase, we build a tripartite hypergraph with the remaining degrees, thus the total of the four types of hyperedges indeed gives a realization of the prescribed degree sequence. The four types of hyperedges, thus the four phases are the following:
\begin{enumerate}
 	\item The role of the first type of hyperedges is to satisfy the degrees of at most two vertices. Then by removing these at most two vertices and adjusting the degrees of the remaining vertices, the number of vertices will be dividable by $3$, all degrees will be at least $\frac{2\left\lfloor\frac{n}{3}\right\rfloor^2}{7}+2*\left\lceil\frac{n}{5}\right\rceil+1$  and at most $\frac{5\left\lfloor\frac{n}{3}\right\rfloor^2}{7}$, and the sum of the degrees is still divisible by $3$.
     \item We then split the vertices into three vertex classes of equal size and similar sums. The role of the second type hyperedge is to ensure that the sum of the degrees in each of the vertex classes are in the same modulo class. It turns out that there is at most one hyperedge of the second type. We adjust the degrees and all degrees will be at least $2\frac{\left\lfloor\frac{n}{3}\right\rfloor^2}{7}+2*\left\lceil\frac{n}{5}\right\rceil$  and at most $\frac{5\left\lfloor\frac{n}{3}\right\rfloor^2}{7}$.
    \item The role of the third type of hyperedges is to provide equal sums in the three vertex classes. Any hyperedge of the third type has three incident vertices from the same vertex class. Therefore it is doable to get the same sums in the three vertex classes because by at this point the sum in the three vertex classes are in the same modulo class by $3$, and each hyperedge decreases the sum of the degrees of its vertex class by  $3$. The remaining degrees sum to the same number in each vertex class and each degree will be between $\frac{2\left\lfloor\frac{n}{3}\right\rfloor^2}{7}$ and $\frac{5\left\lfloor\frac{n}{3}\right\rfloor^2}{7}$.
    \item In the last phase, we build a tripartite hypergraph as a realization of the remaining degrees. This is doable by Theorem~\ref{theo:tripartite-bound}. 
\end{enumerate}

  %  {\bf (TODO: shall we split the parts into separate lemmas?)}

    In the first phase, we consider $(n\mod 3)$ (that is, $0$, $1$ or $2$) arbitrary vertices, call them \emph{extra}. 
    Take $n-2$ vertices that are not extra, and set up a $\left\lceil\frac{n}{5}\right\rceil$-regular degree sequence on it, and subtract $1$ from one of the degrees if both $n-2$ and $\left\lceil\frac{n}{5}\right\rceil$ are odd (thus, the sum of the degrees will be even). This is a regular or almost regular degree sequence that is graphic if $\left\lceil\frac{n}{5}\right\rceil \le n-3 [= (n-2)-1]$ \cite{emt18}, which holds for any $n > 3$. Take a realization of the mentioned degree sequence, $G= (V,E)$. For each extra vertex $v$ extend it into a hypergraph realization $H = (V',E')$, such that $V' = V \cup \{v\}$ and $E' = \{e\cup\{v\}| e\in E\}$. 
    From this hypergraph $H$, we select as many hyperedges as the degree of $v$. This is possible as $(n-2)\left\lceil\frac{n}{5}\right\rceil -1 \ge 2\left\lfloor\frac{5\left\lfloor\frac{n}{3}\right\rfloor^2}{7}\right\rfloor$ for any $n> 2$. For each extra vertex, the number of incident hyperedges of the first type will be its prescribed degree, so after subtracting it, the remaining degree will be 0, and thus, these vertices need no more hyperedges. For other vertices, each of them has at most $\left\lceil\frac{n}{5}\right\rceil$ incident hyperedges for each extra vertex, thus their remaining degree is decreased by at most $2*\left\lceil\frac{n}{5}\right\rceil$. That is, each degree is at least $\frac{2\left\lfloor\frac{n}{3}\right\rfloor^2}{7}+2*\left\lceil\frac{n}{5}\right\rceil+1$ and at most $\frac{5\left\lfloor\frac{n}{3}\right\rfloor^2}{7}$, and the sum of the degrees is still divisible by $3$.

In the second phase, we first split the non-extra vertices into three vertex classes, each containing $\left\lfloor\frac{n}{3}\right\rfloor$ vertices. Let these vertex classes be $A$, $B$, and $C$, let $x_a$, $x_b$, and $x_c$ be the sum of the degrees in these vertex classes, and let $x$ be the sum of all (adjusted) degrees in these vertex classes divided by $3$. We define $d := |x_a-x|+|x_b-x|+|x_c-x|$. While $d$ is larger than $\frac{6\left\lfloor\frac{n}{3}\right\rfloor^2}{7}$, we do the following: Let $d_1$ be the largest degree in the vertex class with the largest degree sum, and let $d_2$ be the smallest degree in the vertex class with the smallest degree sum. It is easy to see that $d_1>d_2$. We swap these two degrees, thus we decrease the largest degree sum and increase the smallest degree sum. We claim that this decreases $d$. Indeed, either there is at most one vertex class whose degree sum is larger than the average, or there is at most one vertex class whose degree sum is smaller than the average. WLOG, assume that there is one vertex class whose degree sum is larger than the average, again WLOG, assume it is $A$ . A moment of thought tells that in this case $d= 2(x_a-x)$. That is $x_a-x> \frac{3\left\lfloor\frac{n}{3}\right\rfloor^2}{7}$. 
The largest decrease in $x_a$ we can make is $\left\lfloor\frac{5\left\lfloor\frac{n}{3}\right\rfloor^2}{7}\right\rfloor - \left \lceil\frac{2\left\lfloor\frac{n}{3}\right\rfloor^2}{7}\right\rceil \ge \frac{3\left\lfloor\frac{n}{3}\right\rfloor^2}{7}$.
Therefore, after the swap, $x_a$ will remain above the average, thus $|x_a-x|$ decreases by the swap. The degree sum in the vertex class of the smallest degree sum is at least half as far from $x$ in absolute value than $x_a$, thus the swap cannot increase the the difference from $x$ in absolute value. What follows is, that $d$ is strictly monotonously decreasing, therefore after finite number of steps, it will be below $\frac{6\left\lfloor\frac{n}{3}\right\rfloor^2}{7}$. Then the difference between the largest degree sum and the smallest degree sum of the vertex classes is also at most $\frac{6\left\lfloor\frac{n}{3}\right\rfloor^2}{7}$.This split of the degrees into three vertex classes is the final split. 

Now either the degree sums are congruent with each other modulo $3$, or they are in three different modulo classes modulo $3$. In the first case, there is no hyperedge of the second type. In the second case, consider the two vertex classes whose degree sums are not minimal. One of them differs from the minimal sum by $1$ modulo $3$, and the other by $-1$ modulo $3$. Choose one and two arbitrary vertices from these two vertex classes, respectively, and add a hyperedge incident with these vertices. By adjusting the degrees, the degree sums will be congruent with each other modulo $3$, and the difference between the maximal and the minimal sum is still at most  $\frac{6\left\lfloor\frac{n}{3}\right\rfloor^2}{7}$. Further each degree is at least $\frac{2\left\lfloor\frac{n}{3}\right\rfloor^2}{7}+2*\left\lceil\frac{n}{5}\right\rceil$  and at most $\frac{5\left\lfloor\frac{n}{3}\right\rfloor^2}{7}$.

In the third phase, we consider those vertex classes whose degree sums are not minimal and consider the differences between their degree sums and the minimal degree sum. These differences are at most $\frac{6\left\lfloor\frac{n}{3}\right\rfloor^2}{7}$, as we discussed.
We spit the $\left\lfloor\frac{n}{3}\right\rfloor$ vertices into $3$ vertex classes, each containing $\left\lfloor\frac{\left\lfloor\frac{n}{3}\right\rfloor}{3}\right\rfloor$ vertices (we might skip at most $2$ vertices). We create a $2\left\lceil\frac{n}{5}\right\rceil$-regular tripartite graph on it,. It is doable for any $n\ge 45$ vertices as in such cases, $2\left\lceil\frac{n}{5}\right\rceil\le \left\lfloor\frac{\left\lfloor\frac{n}{3}\right\rfloor}{3}\right\rfloor^2$, and then we can apply Theorem~\ref{theo:hyper-almost-regular}. Then the degree sum is $2\left\lceil\frac{n}{5}\right\rceil\left\lfloor\frac{\left\lfloor\frac{n}{3}\right\rfloor}{3}\right\rfloor* 3$ which is larger than $\frac{6\left\lfloor\frac{n}{3}\right\rfloor^2}{7}$ if $n\ge17$. We delete the appropriate number of edges from this partite, regular hypergraph such that the remaining degree sum is exactly the prescribed difference. Since the prescribed difference can be divided by $3$, it is doable. Further, each degree is decreased by at most $ 2\left\lceil\frac{n}{5}\right\rceil$. Therefore, now in each vertex class, the degree sum is the same, and each degree is at least $\frac{2\left\lfloor\frac{n}{3}\right\rfloor^2}{7}$ and at most $\frac{5\left\lfloor\frac{n}{3}\right\rfloor^2}{7}$.

Therefore, in the last phase, we can create a tripartite hypergraph, according to Theorem~\ref{theo:tripartite-bound}. 

The hyperedges in the four phases cannot be parallel. Indeed, only the hyperedges in the first phase are incident with extra vertices. Only the hyperedge in the second phase spans exactly two vertex classes. Only the hyperedges in the third phase are inside one vertex class. Finally, only the hyperedges in the fourth phase span all the vertex classes.
Therefore, if we merge the hyperedges introduced in the four phases, we get a realization of the prescribed degree sequence.
\end{proof}

\subsection{Algorithmic considerations}

In this subsection, we explain how the constructive proof of Theorem~\ref{theo:hypergraph-linear-bound} provides a polynomial running time algorithm. Just like in subsection~\ref{ssec:algorithmic-tripartite}, we do not give a thorough analysis but just give reasoning that constructing a hypergraph realization of a degree sequence on $n$ vertices obeying the prescribed constraints can be done in polynomial time.

Constructing the hyperedges in the first phase needs the construction of a realization of an almost regular degree sequence, which can be done in polynomial time with the Havel-Hakimi algorithm. Splitting the vertices into almost equal sums can also be done in polynomial time. Indeed, we only need to compute some sums, find the maximal and minimal degrees, and swap vertices between vertex classes. In each step, a monovariant $d$ is strictly decreasing and it has a starting value bounded above by a polynomial. Furthermore, it can only take integer values and cannot be smaller than $0$. As a result, the number of iterations is bounded above by a polynomial function of $n$. As each iteration can be done in polynomial time, the entire procedure can be done in polynomial time.

Constructing the hyperedge in the second phase needs some easy computations on modulo $3$ and can be done in polynomial time.

Constructing the hyperedges in the third phase needs the construction of a tripartite regular hypergraph, which again can be done in polynomial time as we have discussed in subsection~\ref{ssec:algorithmic-tripartite}.

Then in the last phase, we need to construct another tripartite hypergraph, which can be done in polynomial time as we have discussed in subsection~\ref{ssec:algorithmic-tripartite}.

\subsection{Discussion on the limits}
It is easy to prove the following observation based on the Erd{\H o}s-Gallai theorem \cite{ErdosGallai} and Theorem~\ref{theo:balancing}.
\begin{observation}\label{obs:simple-bound}
    For any $\varepsilon > 0$, there exists an $n_0$ such that for any $n\ge n_0$, any degree sequence on $n$ vertices with even degree sum and each degree between $\left(\frac{1}{4}+\varepsilon\right)n$ and $\left(\frac{3}{4}-\varepsilon\right)n$ is graphic.
\end{observation}
The difference between Observation~\ref{obs:simple-bound} and Theorem~\ref{theo:hypergraph-linear-bound} is much larger than the difference between Observation~\ref{obs:bipartite-bound} and Conjecture~\ref{con:tripartite-bound} (or even between  Observation~\ref{obs:bipartite-bound} and Theorem~\ref{theo:tripartite-bound}). Even more, the bounds in Theorem~\ref{theo:hypergraph-linear-bound} is not in the form $[c\frac{n^2}{2},(1-c)\frac{n^2}{2}]$ (recall that the maximum degree in a simple graph is $n-1$ while it is ${n-1\choose 2}$ in a $3$-uniform hypergraph). Therefore, it is natural to conjecture that the gap in the bounds proved in Theorem~\ref{theo:hypergraph-linear-bound} is far from the largest possible gap. However, it is unclear what the true limit could be. Therefore, we set up a weaker conjecture here.
\begin{conjecture}
    There exists an $\varepsilon>0$ and $n_0\in \mathbb{N}$, such that any hypergraph degree sequence on $n\ge n_0$ vertices, with each degree between $\left(\frac{1}{4}-\varepsilon\right)n^2$ and $\left(\frac{1}{4}+\varepsilon\right)n^2$ and the degree sum divisible by $3$ is graphic.
\end{conjecture}

\section*{Acknowledgments}
IM was supported by NKFIH grant K132696. The project is a continuation of the work done at the Budapest Semesters in Mathematics in 2023 summer. Both authors thank Mufeng Liu and Theodore Mollano for participating in the very early phase of the research at the BSM leading to the results presented in this paper. Thanks to G\'abor Lippner and Tam\'as Keleti for fruitful discussions. Both authors thank BSM for running the program. Both authors thank Dr. Maribel Bueno and Geneva Schlafly for suggesting RL participate in the BSM program and its research opportunities.

\nolinenumbers

% Either type in your references using
% \begin{thebibliography}{}
% Text
% \end{thebibliography}
%
% or
%
% Compile your BiBTeX database using our plos2015.bst
% style file and paste the contents of your .bbl file
% here. See http://journals.plos.org/plosone/s/latex for 

\bibliography{twoseventh}{}
\bibliographystyle{natbib}

\end{document}